\documentclass{amsart}
 
\usepackage{amssymb,latexsym,amsmath,amsfonts,hhline,latexsym,mathrsfs}
\usepackage{pdfsync,verbatimbox}

\newcommand{\cal}{\mathcal}
\usepackage{cite}
\usepackage[matrix,arrow,curve]{xy}
\input xy
\xyoption{all}

\oddsidemargin=\evensidemargin \advance\topmargin by-\headheight

\DeclareMathOperator{\aut}{Aut}

\DeclareMathOperator{\GL}{GL}

\DeclareMathOperator{\orb}{Orb}

\DeclareMathOperator{\reg}{Reg}

\DeclareMathOperator{\Span}{span}

\DeclareMathOperator{\sym}{Sym}

\def\mF{{\mathbb F}}

\def\frC{{\mathscr C}}

\def\cE{{\cal E}}

\def\lg{\langle}

\def\rg{\rangle}

\def\proof{{\bf Proof}.\ }
\def\eprf{\hfill$\square$}
\def\qaq{\quad\text{and}\quad}
\newcommand{\grp}[1]{\langle {#1}\rangle}

\makeatletter
\renewcommand{\subsection}{\@startsection{subsection}{2}{0mm}{-2mm}{-2mm}{\bf\normalsize}}

\makeatother

\newtheorem{formula}{}[section]
\newtheorem{proposition}[formula]{Proposition}
\newtheorem{definition}[formula]{Definition}
\newtheorem{corollary}[formula]{Corollary}
\newtheorem{remark}[formula]{Remark}
\newtheorem{lemma}[formula]{Lemma}
\newtheorem{theorem}[formula]{Theorem}

\newtheorem{example}[formula]{Example}

\def\thrm{\begin{theorem}}
\def\thrml#1{\begin{theorem}\label{#1}}
\def\ethrm{\end{theorem}}
\def\prpstn{\begin{proposition}}
\def\prpstnl#1{\begin{proposition}\label{#1}}
\def\eprpstn{\end{proposition}}
\def\rmrk{\begin{remark}}
\def\rmrkl#1{\begin{remark}\label{#1}}
\def\ermrk{\end{remark}}
\def\dfntn{\begin{definition}}
\def\dfntnl#1{\begin{definition}\label{#1}}
\def\edfntn{\end{definition}}
\def\nmrt{\begin{enumerate}}
\def\enmrt{\end{enumerate}}
\def\tm#1{\item[{\rm (#1)}]}
\def\qtn{\begin{equation}}
\def\qtnl#1{\begin{equation}\label{#1}}
\def\eqtn{\end{equation}}
\def\lmm{\begin{lemma}}
\def\lmml#1{\begin{lemma}\label{#1}}
\def\elmm{\end{lemma}}
\def\crllr{\begin{corollary}}
\def\crllrl#1{\begin{corollary}\label{#1}}
\def\ecrllr{\end{corollary}}
\def\hpthss{\begin{hypothesis}}
\def\hpthssl#1{\begin{hypothesis}\label{#1}}
\def\ehpthss{\end{hypotxesis}}
\def\xmpl{\begin{example}}
\def\xmpll#1{\begin{example}\label{#1}}
\def\exmpl{\end{example}}
\def\css{\begin{cases}}
\def\ecss{\end{cases}}

\begin{document}

\title[A family of permutation groups]{A family of permutation groups with exponentially many
non-conjugated regular elementary abelian subgroups}
\author{\fbox{Sergei Evdokimov}}
\address{Steklov Institute of Mathematics at St. Petersburg, Russia}
\email{evdokim@pdmi.ras.ru}
\author{Mikhail Muzychuk}
\address{Netanya Academic College, Netanya, Israel}
\email{muzy@netanya.ac.il}
\author{Ilia Ponomarenko}
\address{Steklov Institute of Mathematics at St. Petersburg, Russia}
\email{inp@pdmi.ras.ru}
\date{}
	
\begin{abstract}
Given a prime $p$, we construct a permutation group containing at least $p^{p-2}$
non-conjugated regular elementary abelian subgroups of order~$p^3$. This gives 
the first example of a permutation group with exponentially many non-conjugated 
regular subgroups.
\end{abstract}

\maketitle

\section{Introduction}
Regular subgroups of permutation groups arise in many natural contexts, for
example, in group factorizations~\cite{LPS}, Schur rings~\cite{Wie1964}, 
Cayley graphs~\cite{B77}, etc. In the present paper, given a group $H$ and 
a permutation group $\Gamma$, we are interested in the number of $\Gamma$-orbits
\qtnl{280616a}
b_H(\Gamma):=|\orb(\Gamma,\reg(\Gamma,H))|
\eqtn
in the action of~$\Gamma$ by conjugation on the set $\reg(\Gamma,H)$ 
of all its regular subgroups isomorphic to $H$.
Using terminology and arguments of~\cite{B77}, one can see that if $\Gamma$ is the automorphism group of an object of a concrete category~$\frC$, then $b_H(\Gamma)$ equals the number of pairwise non-equivalent representations of this object as a Cayley object over~$H$ in~$\frC$. Note that  $H$ is a CI-group with respect to the category $\frC$ if and only if $b_H(\Gamma)=1$ for every group $\Gamma=\aut(X)$, where $X$ is a Cayley object over~$H$ in~$\frC$. As $\frC$, one can take, for example, the category of finite graphs or other combinatorial structures.\medskip

Let $H$ be a cyclic group. Then, obviously, $b_H(\Gamma)\le c(\Gamma)$,
where $c(\Gamma)$ is the number of the conjugacy classes of full cycles contained
in~$\Gamma$. It was proved in~\cite{Mu99} that the latter
number does not exceed $n=|H|$.\footnote{More exactly, under the Classification of Finite Simple Groups, $c(\Gamma)\le\varphi(n)$ where $\varphi$ is the Euler function, ibid.} Thus, in this case
$b_H(\Gamma)\le n$.\medskip 

The simplest non-cyclic case appears when $H$ is an elementary abelian group~$E_{p^2}$. Here, $b_H(\Gamma)\le b_H(P)$ by the Sylow theorem,  where $P$ is a Sylow $p$-subgroup of the group $\Gamma$. To estimate $b_H(P)$, without loss of
generality one can assume that $P$ is a transitive $p$-group of degree $p^2$,
the action of which on some imprimitivity system induces a regular (cyclic) group of order~$p$, i.e., $P$ belongs to the class defined in the same way as the class $\cE_p$ in Theorem~\ref{251015b} below with~$p^3$ replaced by~$p^2$. With the help of the technique developed in Section~\ref{110716a}, one can describe the set $\reg(\Gamma,H)$ (cf. Theorem~\ref{130516a} and Lemma~\ref{230516a}). Then applying \cite[Theorem~6.1]{EP03}, one can prove that $b_H(P)\le p$. Thus, in this case $b_H(\Gamma)\le n$ too.
\medskip

In the above two cases, the number $b_H(\Gamma)$ does not exceed $n$
for all~$\Gamma$. The main result of the present paper (Theorem~\ref{251015b})
shows that in the general case, neither this bound nor even substantially weaker
bounds are valid.

\thrml{251015b}
Let $H=E_{p^3}$, where $p$ is a prime.
Denote by $\cE_p$ the class of all transitive $p$-groups of degree $p^3$,
the action of which on some imprimitivity system 
induces a regular group isomorphic to~$E_{p^2}$. Then there exists
a group $\Gamma\in\cE_p$ such that $b_H(\Gamma)\ge p^{p-2}$. 
\ethrm

From Theorem~\ref{251015b}, it follows that there is no function $f$, for which the inequality $b_H(\Gamma)\le n^{f(r)}$ holds for all abelian groups $H$ of rank at most~$r$ and all permutation groups~$\Gamma$ of degree~$n$. It would be interesting to find an invariant $t=t(\Gamma)$ such that
$$
b_H(\Gamma)\le n^{f(r,t)}
$$ 
for a function~$f$ in $r$ and $t$; for instance, one can try to take $t=t(\Gamma)$ to be the minimal positive integer $t'$, for which the group~$\Gamma$ is $t'$-closed as a permutation group in the sense of~\cite{W1969}.\footnote{Here for groups $\Gamma\in\cE_p$, the upper bound in inequality~\eqref{090616a} could be useful.}\medskip

The proof of Theorem~\ref{251015b} is given in Section~\ref{140716a}. It is based on a representation of the groups belonging to the class $\cE_p$ with the help of two-variable polynomials over the field~$\mF_p$. The details are presented in~Section~\ref{110716a}. It is interesting to note that for every group~$\Gamma$ from Theorem~\ref{251015b}, the stabilizer of the imprimitivity system is, up to language, a Generalized Reed-Muller code~\cite{KLP}.\medskip

{\bf Notation.}
As usual, $\mF_p$ and $\sym(V)$ denote the field of order~$p$ and the symmetric group on the set~$V$. An elementary abelian $p$-group of order $p^n$ is denoted by~$E_{p^n}$.

\section{Permutation groups and polynomials}\label{110716a}
Let $p$ be a prime. Denote by $R_p$ the factor ring of the polynomial
ring $\mF_p[X,Y]$ modulo the ideal generated by polynomials
$X^p-1$ and $Y^p-1$. The images of the variables $X$ and $Y$ are denoted by
$x$ and $y$, respectively. Denote by $V$ the disjoint union of 
one-dimensional subspaces
$$
V_{i,j}=\{\alpha x^iy^j:\ \alpha\in\mF_p\},\quad i,j=0,\ldots,p-1,
$$
of the ring $R_p$ considered as a linear space over $\mF_p$.\medskip

Every element $f=\sum_{i,j}\alpha_{i,j}x^iy^j$ of $R_p$ yields a permutation 
$$
\sigma_f:\alpha x^iy^j\mapsto (\alpha+\alpha_{i,j})x^iy^j
$$
of the set~$V$. This produces a permutation group on $V$ with $p^2$ orbits $V_{i,j}$
that is isomorphic to the additive group of the ring~$R_p$. For a subgroup $I$ of
the latter group, the corresponding subgroup of $\sym(V)$  is denoted by
$\Delta(I)$. In addition, we define two commuting permutations
$$
\tau_x:\alpha x^iy^j\mapsto\alpha x^{i+1}y^j,\qquad
\tau_y:\alpha x^iy^j\mapsto\alpha x^iy^{j+1}.
$$
Clearly, each of them commutes with the permutation $s=\sigma_{f_0}$, where $f_0=\sum_{i,j}x^iy^j$. The following statement is straightforward.

\lmml{140616a}
In the above notation, we have
\nmrt
\tm{1} $\tau_x^{-1}\sigma_f\tau_x^{}=\sigma_{fx}$ and 
$\tau_y^{-1}\sigma_f\tau_y^{}=\sigma_{fy}$ for all $f\in R_p$,
\tm{2} $G_0:=\grp{s,\tau_x,\tau_y}$ is a regular group on $V$ isomorphic to $E_{p^3}$.  
\enmrt
\elmm

Set $\Gamma(I)$ to be the group generated by $\Delta(I)$ and $\tau_x,\tau_y$.
If $I$ is an ideal of $R_p$, then in view of statement~(1) of Lemma~\ref{140616a}, 
\qtnl{170616c}
\Delta(I)\trianglelefteq \Gamma(I)\qaq\Gamma(I)/\Delta(I)\cong E_{p^2}.
\eqtn
If $I$ is not an ideal, then $\Gamma(I)=\Gamma(I')$, where $I'$ is the ideal of $R_p$ generated by~$I$.

\thrml{130516a}
Let $p$ be a prime. Then
\nmrt
\tm{1} for every ideal $I\ne 0$ of the ring $R_p$, the group $\Gamma(I)$ belongs
to the class $\cE_p$,
\tm{2} every group $\Gamma\in\cE_p$ with $b_H(\Gamma)>0$ is permutation
isomorphic to the group $\Gamma(I)$ for some ideal $I$ of $R_p$.
\enmrt
\ethrm
\proof To prove statement~(1), let $I\ne 0$ be an ideal of $R_p$. Then at least
one of the sets $V_{i,j}$ is an orbit of the group $\Delta(I)$. Since 
$\tau_x$ and $\tau_y$ commute, the group $\grp{\tau_x,\tau_y}$ 
acts regularly on the set~$S=\{V_{i,j}:\ i,j=0,\ldots,p-1\}$. This implies
that the group $\Gamma(I)$ is transitive and $S$ is an imprimitivity 
system of it.
The action of $\Gamma(I)$ on this system induces a regular group isomorphic
to $E_{p^2}$ that is generated  by the images of $\tau_x$ and $\tau_y$
with respect to this action. Thus, $\Gamma(I)\in\cE_p$.\medskip

Let $\Gamma\in\cE_p$. Then $\Gamma$ is a transitive $p$-group of degree $p^3$,
the action of which on some imprimitivity system $S'$ induces a regular group 
isomorphic to~$E_{p^2}$. Without loss of generality, we may assume that
$\Gamma\le\sym(V)$ with $V$ as above. Furthermore, since $b_H(\Gamma)>0$,
the group $\Gamma$ contains a regular subgroup $G'$ isomorphic to $H=E_{p^3}$.
Choose an element $s'\in G'$ such that $\orb(\grp{s'},V)=S'$. Then there
exists a group isomorphism
$$
\varphi:G'\to G_0
$$
taking $s'$ to $s$ (see statement~(2) of Lemma~\ref{140616a}). Since
$\varphi$ is induced by a permutation of~$V$, we may assume
that $S'=S$ and $G_0\in\reg(\Gamma,E_{p^3})$. Note that the permutation $s$
belongs to the stabilizer $\Delta$ of the blocks $V_{i,j}$ in $\Gamma$.
Therefore, $\orb(\Delta,V)=S$. Since the restriction of~$\Delta$ to $V_{i,j}$ is a $p$-group of degree $p$ that contains the restriction of~$s$ to~$V_{i,j}$ for all~$i,j$, this implies that
$$
\Delta\le \Delta(R_p).
$$
It follows that $\Delta=\Delta(I)$ for a subgroup $I$ of $R_p$. Taking
into account that $\Delta$ is normalized by $\tau_x$ and $\tau_y$, we
conclude that $I$ is an ideal of $R_p$ by statement~(1) of 
Lemma~\ref{140616a}.\eprf\medskip

Any maximal element in the class  $\cE_p$  is permutation isomorphic to
the (imprimitive) wreath product of regular groups isomorphic to $E_p$ 
and $E_{p^2}$. One of these maximal elements equals the group 
$\Gamma_p:=\Gamma(R_p)$; set also $\Delta_p=\Delta(R_p)$.
We need two auxiliary lemmas.

\lmml{260516a}
Let $g,h\in R_p$. Then the order of the permutation $t_{g,x}=\sigma_g\tau_x$
(resp. $t_{h,y}=\sigma_h\tau_y$) equals $p$ if and only if $g\in aR_p$
(resp. $h\in bR_p$), where  $a=x-1$ and $b=y-1$.
\elmm
\proof Let $g=\sum_{i,j}\alpha_{i,j}x^iy^j$, and let $v=\alpha x^iv^j$ be a point
of $V$. Then  by the definition of $t_{g,x}$, we have
$$
v^{t_{g,x}}=(\alpha+\alpha_{i,j}) x^{i+1}y^j.
$$
This implies that the order of $t_{g,x}$ equals $p$ if and only if the
following condition is satisfied:
\qtnl{160616a}
\sum_{i=0}^{p-1}\alpha_{i,j}=0,\qquad j=0,\ldots, p-1.
\eqtn
Note that this  is always true, whenever $g\in aR_p$. Conversely,
suppose that relations~\eqref{160616a} hold for some $g\in R_p$. Then 
$$
\alpha_{0,j}=\alpha'_{1,j}-\alpha'_{0,j},\ \ldots,\ 
\alpha_{p-1,j}=\alpha'_{0,j}-\alpha'_{p-1,j},
$$
where  $\alpha'_{i,j}=\sum_{k=0}^{i-1}\alpha_{k,j}$ for all $i,j$. It follows 
that $g=ag'$ with $g'=\sum_{i,j}\alpha'_{i,j}x^iy^j$. This completes
the proof of the first statement. The second statement (on the order of $t_{h,y}$)
is proved similarly.\eprf\medskip

\lmml{230516a}
A permutation group $G$ belongs to the set $\reg(\Gamma_p,E_{p^3})$ if and only if there exist elements $g\in aR_p$ and $h\in bR_p$ such that
\qtnl{230516c}
G=\lg s, t_{g,x},t_{h,y}\rg \qaq ah=bg.
\eqtn
\elmm
\proof To prove the ``only if'' part, suppose that $G\in\reg(\Gamma_p,E_{p^3})$.
Then  $G$ is a self-centralizing subgroup of $\sym(V)$. On the other hand,
the centralizer of $G$ in $\sym(V)$ contains the central element $s$ of the group $\Gamma_p$. Thus, $s\in G$. The other two generators of $G$ can obviously be
chosen so that their images with respect to the epimorphism 
$\Gamma_p\to \Gamma_p/\Delta_p$ coincide with $x$ and $y$. By
Lemma~\ref{260516a}, this implies that there
exist  $g\in aR_p$ and $h\in bR_p$, for which the first equality
in~\eqref{230516c}
holds. Next, since the group $G$ is abelian, the definition of $t_{g,x}$ and
$t_{h,y}$ implies that
$$
\sigma_g\tau_x\,\sigma_h\tau_y=t_{g,x}t_{h,y}=
t_{h,y}t_{g,x}=\sigma_h\tau_y\,\sigma_g\tau_x.
$$
Each of the permutations on the left- and right-hand sides takes
the point $\alpha x^iy^j\in V_{i,j}$ to a certain point 
$\alpha' x^{i+1}y^{j+1}\in V_{i+1,j+1}$. Calculating the images of the former
point with respect to them, we obtain
$$
\alpha+g_{i,j}+h_{i+1,j}=\alpha'=\alpha+h_{i,j}+g_{i,j+1}
$$
or,
equivalently, $h_{i+1,j}-h_{i,j}=g_{i,j+1}-g_{i,j}$ for all $i,j$. Therefore, 
$ah=xh-h=yg-g=bg$, as required.\medskip

Conversely, let $G$ be the group defined by relations~\eqref{230516c}. Then the 
above argument shows that the permutations $s$, $t_{g,x}$, and $t_{h,y}$ 
pairwise commute.
Therefore, the group $G$ is abelian. Moreover, the definition of~$s$ and 
Lemma~\ref{260516a}   imply
that $G$ is elementary abelian and transitive. Thus,
$G\in\reg(\Gamma_p,E_{p^3})$, as required.\eprf

\section{Proof of Theorem~\ref{251015b}}\label{140716a}

By statement~(1) of Theorem~\ref{130516a}, we may restrict ourselves to looking
for the group $\Gamma$ of the form $\Gamma(I)$, where $I$ is an ideal of the 
ring $R_p$.\medskip

For every integer~$k\ge 0$, set
$$
I_k=\Span_{\mF_p}\{a^ib^j:\ i+j\ge k\},
$$
where the elements $a$ and $b$ are as in Lemma~\ref{230516a}. Clearly, $I_k$ is an
ideal of $R_p$, and $I_{k+1}\subseteq I_k$ for all~$k$, and also $I_k=0$ for $k> 2(p-1)$. Below, the kernels of the mappings
$I_k\to aI_k$ and $I_k\to bI_k$ induced by the multiplication
by $a$ and $b$ are denoted by $A_k$ and $B_k$, respectively.

\lmml{240516f}
Suppose that $p\le k\le 2(p-1)$. Then
\nmrt
\tm{1}  $\dim(I_k)={{2p-k}\choose{2}}$, 
\tm{2}  $aI_k=bI_k=I_{k+1}$,
\tm{3}  $\dim(A_k)=\dim(B_k)=2p-k-1$.
\enmrt
\elmm
\proof The leading monomials of the polynomials 
$$
(x-1)^i(y-1)^j,\qquad 0\le i,j\le p-1,
$$ 
are obviously 
linearly independent. Therefore, the polynomials $a^ib^j$ with $i+j\ge k$
form a linear basis of the ideal $I_k$. This immediately proves statement~(1).
To prove statement~(2), we note that, obviously, $aI_k\subseteq I_{k+1}$.
Conversely, let $c\in I_{k+1}$. Since $k\ge p$, we have
$c=abu$ for some $u\in I_{k-1}$, which proves the reverse
inclusion. The rest of statement~(2) is proved similarly.
Finally, statement~(3) follows, because the linear space $A_k$ (resp.~$B_k$)
is spanned by the monomials $a^{p-1}b^i$ (resp. $a^ib^{p-1}$) with 
$k-p+1\le i\le p-1$.\eprf\medskip

In what follows, for a subgroup $G$ of a group $\Gamma$ we denote by $G^\Gamma$ the set of all $\Gamma$-conjugates of~$G$.

\lmml{070616a}
Let $\Gamma_{k,p}=\Gamma(I_k)$,  where $k$ is as in Lemma~\ref{240516f}. Then
\nmrt
\tm{1} $|\Gamma_{k,p}|=p^{2+\dim(I_k)}$,
\tm{2} $|\reg(\Gamma_{k,p},E_{p^3})|=p^{\dim(A_k)+\dim(B_k)+\dim(I_{k+1})-2}$,
\tm{3} $p^{\dim(I_k)-4}\le |G^{\Gamma_{k,p}}|\le p^{\dim(I_k)-1}$ for all $G\in\reg(\Gamma_{k,p},E_{p^3})$.
\enmrt
\elmm
\proof Obviously, $|\Delta(I_k)|=p^{\dim(I_k)}$. Therefore, statement~(1)
follows from the right-hand side of formula~\eqref{170616c}. Next, from
Lemma~\ref{230516a} it follows that
$$
\reg(\Gamma_{k,p},E_{p^3})=\{G_{g,h}:\ (g,h)\in M\},
$$
where $G_{g,h}=\grp{s,t_{g,x}, t_{h,y}}$ and
\qtnl{170616d}
M=\{(g,h)\in (I_k\cap aR_p)\times (I_k\cap bR_p):\ ah=bg\}.
\eqtn
However, $I_k\cap aR_p=I_k\cap bR_p=I_k$, because $k\ge p$. So 
by statement~(2) of Lemma~\ref{240516f}, the element $ah=bg$ runs over
the ideal $I_{k+1}$, whenever $(g,h)$ runs over the set~$M$.  By the
definition of $A_k$ and $B_k$, this implies that
$$
|M|=p^{\dim(A_k)+\dim(B_k)+\dim(I_{k+1})}.
$$
Thus to complete the proof of statement~(2), it suffices to verify that
$G_{g^{},h^{}}=G_{g',h'}$ if and only if $t_{g^{},x^{}}=s^it_{g',x}$ and
$t_{h^{},y^{}}=s^jt_{h',y^{}}$ for some $0\le i,j\le p-1$. However, this is true,
because $G_{g^{},h^{}}=G_{g',h'}$ if and only if
$\varphi(G_{g^{},h^{}})=\varphi(G_{g',h'})$, where $\varphi$ is the quotient
epimorphism of~$\Gamma_{k,p}$ modulo the group~$\grp{s}$.\medskip

To prove statement~(3), we note that in view of statement~(1),
\qtnl{230616a}
|G^\Gamma|=\frac{|\Gamma|}{|N|}=
\frac{p^{2+\dim(I_k)}}{|C|\cdot|N/C|},
\eqtn
where $\Gamma=\Gamma_{k,p}$, and $N$ and $C$ are, respectively, the normalizer and centralizer of~$G$ in~$\Gamma$. Since $G$ is a regular elementary abelian group and the quotient $N/C$ is isomorphic to a subgroup of  a Sylow $p$-subgroup $P$ of the group $\aut(G)\cong\GL(3,p)$ (here we use the fact that $\Gamma$ is a $p$-group), we conclude that
$$
|C|=|G|=p^3\qaq 
1\le |N/C|\le |P|.
$$ 
However, $|P|=p^3$. Thus, statement~(3) follows  from
formula~\eqref{230616a}.\eprf\medskip

To complete the proof of Theorem~\ref{251015b}, we note that
$\reg(\Gamma_k,E_{p^3})$ is the disjoint union of 
distinct sets $G^{\Gamma_k}$,  where $\Gamma_k=\Gamma_{k,p}$ as in Lemma~\ref{070616a} and $G\in\reg(\Gamma_k,E_{p^3})$. 
Therefore, setting  $m_k$ and $M_k$ to be, respectively, the minimum and maximum 
of the numbers $|G^{\Gamma_k}|$, we obtain
\qtnl{100616a}
\frac{|\reg(\Gamma_k,E_{p^3})|}{m_k}\ge 
b_H(\Gamma_k)\ge \frac{|\reg(\Gamma_k,E_{p^3})|}{M_k}.
\eqtn
However, by statement~(3) of Lemma~\ref{070616a}, 
$
m_k\ge p^{\dim(I_k)-4}$ and $M_k\le p^{\dim(I_k)-1}$.
By statement~(2)  of Lemma~\ref{070616a}, this implies that
$$
\frac{|\reg(\Gamma_k,E_{p^3})|}{m_k}\le 
p^{d+2}
\qaq
\frac{|\reg(\Gamma_k,E_{p^3})|}{M_k}\ge 
p^{d-1}.
$$
where $d=\dim(A_k)+\dim(B_k)+\dim(I_{k+1})-\dim{I_k}$.
Besides, by statements~(1) and~(3) of Lemma~\ref{240516f}, we have
$d=2p-k-1$.
Thus,  
\qtnl{090616a}
p^{2p-k+1}\ge b_H(\Gamma_k)\ge p^{2p-k-2}.
\eqtn
This lower bound for $b_H(\Gamma_k)$ with $k=p-1$ proves Theorem~\ref{251015b}.

\end{document}